\documentclass[12pt,a4paper]{article}
\usepackage[margin=27mm]{geometry}
\usepackage[T1]{fontenc}
\usepackage{lmodern,amsmath,amssymb,amsthm,microtype}
\usepackage{needspace}
\usepackage{xcolor}
\usepackage[hidelinks]{hyperref}
\newtheorem{theorem}{Theorem}[section]
\newtheorem{lemma}[theorem]{Lemma}
\newtheorem{proposition}[theorem]{Proposition}
\newtheorem{corollary}[theorem]{Corollary}
\theoremstyle{definition}
\newtheorem{definition}[theorem]{Definition}
\theoremstyle{remark}
\newtheorem{remark}[theorem]{Remark}
\DeclareMathOperator{\mis}{mis}
\DeclareMathOperator{\itm}{itm}
\newcommand{\M}{\mathcal M}
\newcommand{\B}{\mathcal B}
\newcommand{\ind}{\mathbf 1}
\title{Maximal independent sets in graphs with\\ bounded minimum degree}
\author{Guorong Gao\thanks{E-mail: \href{mailto:grgao@fzu.edu.cn}{grgao@fzu.edu.cn}.},
Ziyan Zhang\thanks{E-mail: \href{mailto:ziyanzhang1105@163.com}{ziyanzhang1105@163.com}.},
Hongbin Zhao\thanks{E-mail: \href{mailto:hbzhao2024@163.com}{hbzhao2024@163.com}.}\\[0.5em]
\small School of Mathematics and Statistics, Fuzhou University,\\
\small Fuzhou, Fujian, China}
\date{}
\hypersetup{
pdftitle={Maximal independent sets in graphs with bounded minimum degree},
pdfauthor={Guorong Gao, Ziyan Zhang, and Hongbin Zhao}
}
\begin{document}
\maketitle
\begin{abstract}
For integers $k\ge0$ and $n\ge k+1$, let $F(n,k)$ be the maximum number of maximal independent sets in an $n$-vertex graph of minimum degree at least $k$. For every fixed $k\ge3$, we determine $F(n,k)$ for all sufficiently large $n$ and characterize all extremal graphs. The proof combines a uniform counting gap at minimum degree two with a bound in terms of induced triangle matchings and induction by vertex deletion. A finite optimization then determines the additive terms and the equality cases.
\end{abstract}

\section{Introduction}\label{sec:intro}
An independent set is a set of pairwise nonadjacent vertices. A maximal independent set, abbreviated MIS, is an independent set not properly contained in another independent set. This is a condition on inclusion, rather than on cardinality. We write $\M(G)$ for the family of MISs of $G$ and $\mis(G)=|\M(G)|$. All graphs considered here are finite and simple.

The problem of maximizing $\mis(G)$ over graphs of a given order goes back to Miller and Muller~\cite{MillerMuller} and Moon and Moser~\cite{MoonMoser}. The former obtained the count in 1960; the latter studied the equivalent question for maximal cliques, which become MISs upon taking complements. Their result gives $g(0)=g(1)=1$ and, for $n\ge2$,
\begin{equation*}
g(n)=\begin{cases}
3^{n/3},&n\equiv0\pmod3,\\
4\cdot3^{(n-4)/3},&n\equiv1\pmod3,\\
2\cdot3^{(n-2)/3},&n\equiv2\pmod3.
\end{cases}
\end{equation*}
For $n\ge2$, the extremal graphs are disjoint unions of triangles and a bounded number of cliques of orders two or four. Each MIS chooses one vertex from every component. Wood~\cite{Wood} gives a short proof of the bound.

\begingroup
\emergencystretch=1em
\microtypesetup{protrusion=false}
Subsequent work studied how structural restrictions change this maximum. F\"uredi~\cite{Furedi} treated connected graphs of sufficiently large order, and Griggs, Grinstead, and Guichard~\cite{GGG} determined the connected maximum and all equality cases. Hujter and Tuza~\cite{HT} solved the triangle-free problem, while Jin and Li~\cite{JL} determined the second-largest count in the unrestricted problem. Kahn and Park~\cite{KP} showed that graphs whose MIS counts are close to the Moon--Moser bound on the exponential scale contain an induced triangle matching covering nearly all vertices. Palmer and Patk\'os~\cite{PP} strengthened this stability result by determining the exact maximum under a bound on the size of an induced triangle matching. For further related results,  we refer the reader to \cite{BGLY,BLST,GHZ,GDTD,ZGW,STW}.\par
\endgroup

We study the maximum number of MISs under a lower bound on the minimum degree. For integers $k\ge0$ and $n\ge k+1$, put
\begin{equation*}
F(n,k)=\max\{\mis(G): |V(G)|=n,\ \delta(G)\ge k\}.
\end{equation*}

Xiao, Wang, and Tu~\cite{XWT} determined $F(n,3)$ for every $n\ge4$ and characterized all extremal graphs. For sufficiently large $n$, each extremal graph consists of a vertex joined to every vertex of a disjoint union of triangles, together with at most two detached copies of $K_4$. We extend their result to every fixed minimum-degree lower bound $k\ge4$, determining the exact maximum and all equality cases for sufficiently large $n$.

When $k\ge3$, no triangle can occur as a component, but the independent choices in disjoint triangles are preserved if every triangle vertex is joined to a common set of $k-2$ vertices. The resulting construction retains the exponential rate $3^{n/3}$.

We first describe the family that contains all extremal graphs. For $X\subseteq V(G)$, let $G[X]$ denote the subgraph of $G$ induced by $X$. A vertex of a graph $R$ is \emph{universal} if it is adjacent to every other vertex; the vertex of $K_1$ is universal. Let $U(R)$ be the set of universal vertices.

\begin{definition}\label{def:B}
For $s\ge1$, $q\in\{0,1,2\}$, and $t\ge1$, let $\B_{s,q,t}$ consist of graphs admitting a partition
\[
V(G)=S\mathbin{\dot\cup}C_1\mathbin{\dot\cup}\cdots\mathbin{\dot\cup}C_q
\mathbin{\dot\cup}T_1\mathbin{\dot\cup}\cdots\mathbin{\dot\cup}T_t
\]
such that $|S|=s$, each $C_j$ induces $K_4$, and each $T_i$ induces a triangle. There are no edges between distinct blocks among the $C_j$ and $T_i$. Every triangle vertex is adjacent to all of $S$, and every $C_j$ vertex has at least $s-1$ neighbors in $S$. The graph $R=G[S]$ is arbitrary. We call any partition satisfying these conditions a \emph{defining partition} of $G$. We call $S$ the \emph{support} and $J=G[S\cup C_1\cup\cdots\cup C_q]$ the \emph{kernel} of this partition.
\end{definition}

Every triangle vertex has degree $s+2$, every $K_4$ vertex has degree at least $s+2$, and every support vertex has at least $3t$ neighbors. Thus all members of $\B_{s,q,t}$ have minimum degree at least $s+2$ when $3t\ge s+2$.

For a defining partition, let $N_G(u)$ denote the set of neighbors of $u$ in $G$, and define
\begin{equation*}
a_j(u)=|C_j\setminus N_G(u)|,\qquad P=\{u\in S:a_j(u)>0\text{ for some }j\}.
\end{equation*}
Each $C_j$ vertex has at most one nonneighbor in $S$, so $\sum_{u\in S}a_j(u)\le4$.
Define $\gamma_1=0$, $\gamma_2=1$, $\gamma_3=2$ and $\gamma_s=g(s)$ for $s\ge4$.

\begin{theorem}\label{thm:main}
Fix $k\ge3$ and put $s=k-2$. For all sufficiently large $n$, let $q\in\{0,1,2\}$ satisfy $q\equiv n-s\pmod3$ and put $t=(n-s-4q)/3$. Then
\begin{equation}\label{eq:main}
F(n,k)=\begin{cases}
3^t+g(s),&q=0,\\
4^q3^t+4^q+\gamma_s,&q=1,2.
\end{cases}
\end{equation}
The extremal graphs are precisely those described in Theorem~\ref{thm:equality}, and each has minimum degree exactly $k$.
\end{theorem}

\Needspace{8\baselineskip}
\begin{theorem}\label{thm:equality}
With $s,q,t$ as in Theorem~\ref{thm:main}, an $n$-vertex graph is extremal if and only if it admits a defining partition satisfying the corresponding condition below, where $R=G[S]$.
\begin{enumerate}
\renewcommand{\labelenumi}{(\roman{enumi})}
\item If $q=0$, then $\mis(R)=g(s)$.
\item If $q=1$, then $\sum_{u\in S}a_1(u)=4$ and $\mis(R)-|U(R)\cap P|=\gamma_s$.
\item If $q=2$, every vertex of $C_1\cup C_2$ has the same unique nonneighbor in $S$, say $u$, and $\mis(R)-\ind_{\{u\in U(R)\}}=\gamma_s$.
\end{enumerate}
Here $\ind_E$ is $1$ when $E$ holds and $0$ otherwise.
\end{theorem}

We first establish a uniform counting gap for graphs of minimum degree at least two. For sufficiently large graphs of minimum degree at least three whose MIS count is at least a fixed positive proportion of $3^{n/3}$, a bound in terms of induced triangle matchings yields a vertex of large degree. Induction by vertex deletion then shows that every extremal graph belongs to the family in Definition~\ref{def:B}. Finally, optimizing over the possible kernels determines the additive term and all equality cases.

The rest of this paper is organized as follows.
Section 2 collects the counting tools. Section 3 establishes the minimum-degree-two gap, and Section 4 proves the structural gap. Section 5 proves Theorems 1.2 and 1.3 through kernel optimization.

\section{Counting tools}
For $X\subseteq V(G)$, write $G-X=G[V(G)\setminus X]$. A vertex $x$ has open neighborhood $N_G(x)$, closed neighborhood $N_G[x]=N_G(x)\cup\{x\}$, and degree $d_G(x)=|N_G(x)|$. The minimum and maximum degrees are $\delta(G)$ and $\Delta(G)$. For a nonnegative integer $r$, write $rG$ for the disjoint union of $r$ copies of $G$, with $0G=\varnothing$. We use $\mis(\varnothing)=1$, and empty products are one.

An independent set is maximal exactly when it \emph{dominates} the graph: every vertex outside it has a neighbor in it. An undominated vertex could otherwise be added. This also shows that MISs are chosen independently in different components, giving $\mis(G_1\cup\cdots\cup G_h)=\prod_{i=1}^{h}\mis(G_i)$ for a disjoint union.

\begin{theorem}[Miller--Muller; Moon--Moser~\cite{MillerMuller,MoonMoser}]\label{thm:MM}
Every $n$-vertex graph $G$ satisfies $\mis(G)\le g(n)$. For $n\ge2$, the equality graphs, up to isomorphism, are
\[
\begin{array}{rl}
n=3a: & aK_3,\\
n=3a+1: & K_4\cup(a-1)K_3\ \text{or}\ 2K_2\cup(a-1)K_3,\\
n=3a+2: & K_2\cup aK_3.
\end{array}
\]
For $n=0,1$, the unique graph of that order attains equality.
\end{theorem}
In particular, $g(n)\le3^{n/3}$. Also, $g$ is nondecreasing, since adding an isolated vertex preserves the number of MISs.

\begin{theorem}[Griggs--Grinstead--Guichard~\cite{GGG}]\label{thm:GGG}
Every connected graph of order $r\ge6$ has at most $c(r)$ MISs, where
\begin{equation*}
c(r)=\begin{cases}
2\cdot3^{(r-3)/3}+2^{(r-3)/3},&r\equiv0\pmod3,\\
3^{(r-1)/3}+2^{(r-4)/3},&r\equiv1\pmod3,\\
4\cdot3^{(r-5)/3}+3\cdot2^{(r-8)/3},&r\equiv2\pmod3.
\end{cases}
\end{equation*}
\end{theorem}

\begin{proposition}\label{prop:branch}
For every graph $G$ and every vertex $x\in V(G)$,
\begin{equation*}
\mis(G)\le\mis(G-x)+\mis(G-N_G[x]).
\end{equation*}
Moreover, the MISs of $G$ containing $x$ are in bijection with those of $G-N_G[x]$.
\end{proposition}
\begin{proof}
Partition the MISs of $G$ into those containing $x$ and those avoiding $x$. If $I$ contains $x$, then
$I=\{x\}\cup I'$, where $I'$ is an MIS of $G-N_G[x]$. If $I$ avoids $x$, then $I$ is an MIS of $G-x$; every vertex outside $I$ in $G-x$ is already dominated by $I$ in $G$. This proves the inequality.
\end{proof}

An \emph{induced triangle matching} is an induced subgraph whose components are triangles. Its size is the number of triangles. Let $\itm(G)$ denote the maximum size of an induced triangle matching in $G$.

\begin{theorem}[Palmer--Patk\'os~\cite{PP}]\label{thm:PPexact}
Let $n\ge4$ and $0\le t\le\lfloor n/3\rfloor$ be integers, and put $r=n-3t$. The maximum number of MISs in an $n$-vertex graph containing no induced triangle matching of size $t+1$ is
\[
\begin{cases}
3^t2^{r/2},&r\text{ is even},\\
3^{t-1}2^{(r+3)/2},&r\text{ is odd and }t>0,\\
5\cdot2^{(n-5)/2},&r\text{ is odd and }t=0.
\end{cases}
\]
\end{theorem}

We will use the following consequence of Theorem~\ref{thm:PPexact}.

\begin{proposition}\label{prop:PP}
If $G$ has order $n\ge4$ and $\tau=\itm(G)$, then
\begin{equation*}
\mis(G)\le3^\tau 2^{(n-3\tau )/2}=3^{n/3}\rho^{n-3\tau},\qquad \rho=\frac{\sqrt2}{3^{1/3}}<1.
\end{equation*}
\end{proposition}
\begin{proof}
Put $r=n-3\tau$ and apply Theorem~\ref{thm:PPexact} with $t=\tau$. When $r$ is even, the stated bound follows directly. In the other two cases, the ratios of the bounds in Theorem~\ref{thm:PPexact} to $3^\tau 2^{r/2}$ are $2^{3/2}/3$ and $5/2^{5/2}$, respectively, both less than one. The identity follows by substitution, and $\rho<1$ because $8<9$.
\end{proof}

\section{The minimum-degree-two gap}
For $n\ge3$, define
\begin{equation*}
\Pi(n)=\max\left\{\prod_{i=1}^h b_i:h\ge1,\ b_i\in\mathbb Z,\ b_i\ge3,\ \sum_{i=1}^h b_i=n\right\}.
\end{equation*}
Thus $\Pi(n)$ is the maximum number of MISs in a disjoint union of cliques of orders at least three. Concatenating partitions gives $\Pi(a+b)\ge\Pi(a)\Pi(b)$ for $a,b\ge3$.

\begin{lemma}\label{lem:partition}
We have $\Pi(3)=3$, $\Pi(4)=4$, and $\Pi(5)=5$. For $n\ge6$, write $n=4q+3t$, where $q\in\{0,1,2\}$ and $t\ge0$. Then $\Pi(n)=4^q3^t$, and the unique maximizing multiset consists of $q$ fours and $t$ threes.
\end{lemma}
\begin{proof}
The first three values are immediate. Now let $n\ge6$ and consider a maximizing partition. A part $b\ge6$ can be replaced by $3$ and $b-3$, increasing the product because $3(b-3)>b$. Thus every part of a maximizing partition lies in $\{3,4,5\}$. If a part is $5$, then $n\ge6$ ensures that there is another part. If another part is $3$, replace the parts $5$ and $3$ by two parts of size $4$. If it is $4$, replace the parts $5$ and $4$ by three parts of size $3$. If it is $5$, replace the two parts of size $5$ by parts of sizes $4,3,3$. Each replacement strictly increases the product, since $15<16$, $20<27$, and $25<36$. Thus an optimal partition has only threes and fours. Three fours can be replaced by four threes, increasing $64$ to $81$. Thus a maximizing partition has $q\in\{0,1,2\}$ fours. Since $4\equiv1\pmod3$, the identity $n=4q+3t$ determines $q$ uniquely by $q\equiv n\pmod3$, and then determines $t=(n-4q)/3$. Hence the maximizing multiset is unique and $\Pi(n)=4^q3^t$.
\end{proof}

\begin{lemma}\label{lem:basegap}
Let $n\ge6$ and $n=4q+3t$ as above. Every $n$-vertex graph $G$ with $\delta(G)\ge2$ either satisfies $G\cong qK_4\cup tK_3$ or $\mis(G)\le\frac{15}{16}\Pi(n)$.
\end{lemma}
\begin{proof}
Every component has at least three vertices. At orders three and four, Theorem~\ref{thm:MM} gives bounds of three and four, attained among connected graphs only by $K_3$ and $K_4$. For graphs of order five, Theorem~\ref{thm:MM} gives an upper bound of six, attained only by the disconnected graph $K_2\cup K_3$. A connected five-vertex graph therefore contributes at most five. These bounds are $\Pi(3),\Pi(4),\Pi(5)$.

For a connected component of order $r\ge6$, Theorem~\ref{thm:GGG} and Lemma~\ref{lem:partition} give the following ratios, with $a\ge2$:
\[
\begin{aligned}
\frac{c(3a)}{\Pi(3a)}&=\frac23+\frac13\left(\frac23\right)^{a-1}\le\frac89,\\
\frac{c(3a+1)}{\Pi(3a+1)}&=\frac34+\frac14\left(\frac23\right)^{a-1}\le\frac{11}{12},\\
\frac{c(3a+2)}{\Pi(3a+2)}&=\frac34+\frac3{16}\left(\frac23\right)^{a-2}\le\frac{15}{16}.
\end{aligned}
\]
Each ratio is maximized at $a=2$. Let $r_1,\ldots,r_h$ be the component orders of $G$. If a component has order $r_i\ge6$, then it has at most $\frac{15}{16}\Pi(r_i)$ MISs. Using $\Pi(r_j)$ as the bound for each other component, the multiplicativity of mis and the supermultiplicativity of $\Pi$ give $\mis(G)\le\frac{15}{16}\prod_{i=1}^h\Pi(r_i)\le\frac{15}{16}\Pi(n)$.

Now suppose that every component has order $3$, $4$, or $5$. Each component contributes at most its order. If the partition of $n$ into component orders does not attain $\Pi(n)$, one of the replacements in the proof of Lemma~\ref{lem:partition} applies. For the four replacements, the ratios of the old product to the new product are $15/16$, $20/27$, $25/36$, and $64/81$, respectively, all at most $15/16$. The new product is at most $\Pi(n)$, so the original product, and hence $\mis(G)$, is at most $\frac{15}{16}\Pi(n)$.

If the partition into component orders is optimal, it consists of $q$ fours and $t$ threes. Every component of order three is a triangle because $\delta(G)\ge2$. If $G\not\cong qK_4\cup tK_3$, then some component of order four is not complete. By Theorem~\ref{thm:MM}, that component has at most three MISs, so $\mis(G)\le\frac34\Pi(n)$.
\end{proof}

\begin{corollary}\label{cor:degree2}
We have $F(3,2)=3$, $F(4,2)=4$, and $F(5,2)=5$.
For $n\ge6$,
\begin{equation*}
F(n,2)=\Pi(n)=\begin{cases}
g(n),&n\equiv0,1\pmod3,\\
\dfrac89g(n),&n\equiv2\pmod3.
\end{cases}
\end{equation*}
For $n\ge6$, the unique extremal graph, up to isomorphism, is $qK_4\cup tK_3$, where $n=4q+3t$ as above.
\end{corollary}
\begin{proof}
For $n=3,4$, Theorem~\ref{thm:MM} gives $F(n,2)\le n$, and $K_n$ attains this bound. For $n=5$, Theorem~\ref{thm:MM} gives an upper bound of six, attained only by $K_2\cup K_3$ up to isomorphism. Since this graph has minimum degree one, $F(5,2)\le5$, and $K_5$ attains this bound. For $n\ge6$, the graph $qK_4\cup tK_3$ has minimum degree at least two and attains $\Pi(n)$. Lemma~\ref{lem:basegap} therefore gives $F(n,2)=\Pi(n)$ and uniqueness up to isomorphism. By Lemma~\ref{lem:partition}, $\Pi(n)=g(n)$ when $n\equiv0,1\pmod3$. If $n=3a+2$, then $q=2$ and $t=a-2$, so $\Pi(n)=16\cdot3^{a-2}=\frac89\bigl(2\cdot3^a\bigr)=\frac89g(n)$.
\end{proof}
Theorem~\ref{thm:MM} also gives $F(n,0)=g(n)$ for $n\ge1$ and $F(n,1)=g(n)$ for $n\ge2$.

\section{The structural gap}
\begin{lemma}\label{lem:large}
Fix integers $s,D\ge1$ and a real number $\beta>0$. Every sufficiently large graph $G$ satisfying
\[
\delta(G)\ge s+2,\qquad \mis(G)\ge\beta3^{|V(G)|/3}
\]
has a vertex of degree greater than $D$.
\end{lemma}
\begin{proof}
Put $n=|V(G)|$ and suppose $\Delta(G)\le D$. Take a largest induced triangle matching, with $\tau=\itm(G)$ triangles, and let $W$ be the set of uncovered vertices. Every covered vertex has exactly two neighbors within the matching and hence at least $s$ in $W$. Counting edges between the matching and $W$ gives
\[
3s\tau\le D|W|=D(n-3\tau),\qquad n-3\tau\ge\frac{s}{D+s}n.
\]
Since $0<\rho<1$, Proposition~\ref{prop:PP} yields $\mis(G)\le3^{n/3}\rho^{sn/(D+s)}$. For sufficiently large $n$, the last factor is smaller than $\beta$, contradicting the assumed lower bound. For $0<\beta<1$, it suffices to take $n\ge4$ and $n>\frac{D+s}{s}\frac{\log(1/\beta)}{\log(1/\rho)}$. For $\beta\ge1$, the contradiction is immediate once $n\ge4$.
\end{proof}

\begin{lemma}\label{lem:gap}
For each fixed $s\ge1$, there is an integer $t_0(s)$ such that for every $q\in\{0,1,2\}$ and integer $t\ge t_0(s)$, every graph $G$ of order $s+4q+3t$ with $\delta(G)\ge s+2$ either belongs to $\B_{s,q,t}$ or satisfies
\[
\mis(G)\le4^q3^t-c_s3^t,\qquad c_s=2^{-(s+4)}.
\]
\end{lemma}
\begin{proof}
We use induction on $s$. \emph{Base case $s=1$.} Put $n=4q+3t$ and suppose $\mis(G)>\frac{31}{32}\Pi(n)$. The ratio of this lower bound to $3^{(n+1)/3}$ is the positive constant $(31/32)4^q3^{-(4q+1)/3}$. Lemma~\ref{lem:large}, with $D=10$, gives a vertex $x$ of degree at least eleven when $t$ is sufficiently large. Let $H=G-x$. Then $|V(H)|=n$ and $\delta(H)\ge2$, since deleting $x$ reduces the degree of each remaining vertex by at most one. By Proposition~\ref{prop:branch}, Theorem~\ref{thm:MM}, and the monotonicity of $g$,
$\mis(G)\le\mis(H)+g(n-d_G(x))\le\mis(H)+g(n-11)$.
For sufficiently large $n$, direct substitution gives
\[
\frac{g(n-11)}{\Pi(n)}=\begin{cases}
4/243,&n\equiv0\pmod3,\\
1/54,&n\equiv1\pmod3,\\
1/48,&n\equiv2\pmod3,
\end{cases}\qquad\text{each less than }\frac1{32}.
\]
If $H\not\cong qK_4\cup tK_3$, Lemma~\ref{lem:basegap} would imply $\mis(G)<\left(\frac{15}{16}+\frac1{32}\right)\Pi(n)=\frac{31}{32}\Pi(n)$, contrary to our assumption. Hence $H\cong qK_4\cup tK_3$. Every triangle vertex has degree two in $H$ and must be adjacent to $x$, since $\delta(G)\ge3$. When $s=1$, the definition imposes no restriction on the edges from $x$ to the $K_4$ components. Thus $G\in\B_{1,q,t}$ with support $\{x\}$. Taking the contrapositive, if $G\notin\B_{1,q,t}$, then $\mis(G)\le\frac{31}{32}4^q3^t\le4^q3^t-\frac1{32}3^t$, which is the required gap for $c_1=1/32$.

\smallskip\noindent\emph{Inductive step.}
Let $s\ge2$, assume the assertion for $s-1$, and put $c_s=c_{s-1}/2$. Choose an integer $D$ satisfying $3^{(s+7-D)/3}\le c_{s-1}/2$. Suppose $\mis(G)>4^q3^t-c_s3^t$. The ratio of this lower bound to $3^{|V(G)|/3}$ is the positive constant $(4^q-c_s)3^{-(s+4q)/3}$. Lemma~\ref{lem:large} therefore gives $x$ with $d_G(x)>D$. Put $H=G-x$; then $|V(H)|=(s-1)+4q+3t$, $\delta(H)\ge s+1$. By Theorem~\ref{thm:MM}, the monotonicity of $g$, the bound $4q\le8$, and the choice of $D$, for sufficiently large $t$,
\begin{align*}
    \mis(G-N_G[x])
&\le g(s+4q+3t-1-d_G(x)) \le g(s+4q+3t-1-D)\\
&\le 3^t3^{(s+4q-1-D)/3}\le \frac{c_{s-1}}{2}\,3^t.
\end{align*}
If $H\notin\B_{s-1,q,t}$, the induction hypothesis and Proposition~\ref{prop:branch} give
$\mis(G)\le4^q3^t-c_{s-1}3^t+\frac{c_{s-1}}2\,3^t=4^q3^t-c_s3^t$,
contradicting our assumption. Thus $H\in\B_{s-1,q,t}$.

Take a defining partition of $H$ with support $S'$. Each triangle vertex has degree $s+1$ in $H$ and must be adjacent to $x$. Set $S=S'\cup\{x\}$. Then every triangle vertex is adjacent to all of $S$. A vertex $v\in C_j$ has three neighbors in its own block and no neighbors in any other block $C_\ell$ or $T_i$. Its remaining neighbors in $G$ lie in $S$, so $|N_G(v)\cap S|=d_G(v)-3\ge s-1$. All other partition conditions are inherited from $H$. Hence $G\in\B_{s,q,t}$.

For each $q\in\{0,1,2\}$, choose a threshold for $t$ so that the induction hypothesis, Lemma~\ref{lem:large}, and the estimates above apply. Let $t_0(s)$ be the maximum of these three thresholds. This completes the inductive step. Finally, $c_1=2^{-5}$ and $c_s=c_{s-1}/2$ give $c_s=2^{-(s+4)}$.
\end{proof}

\section{Kernel optimization and the main theorem}
By Lemma~\ref{lem:gap}, it remains to optimize over $\B_{s,q,t}$. For fixed $s,q,t$, the variable part of the construction is the kernel, which has $s+4q\le s+8$ vertices. The MISs avoiding the support contribute $4^q3^t$, while those meeting the support contribute a term that depends only on the kernel. To optimize this term, we first give a characterization of $\gamma_s$.

\begin{lemma}\label{lem:gamma}
For $s\ge1$,  $\gamma_s=\max\limits_{\substack{|V(R)|=s\\u\in V(R)}} \left(\mis(R)-\ind_{\{u\in U(R)\}}\right).$
\end{lemma}
\begin{proof}
For $1\le s\le3$, the graph $K_s$ gives $s-1$, since it has $s$ MISs and every vertex is universal. By Theorem~\ref{thm:MM}, every noncomplete graph of these orders has at most $s-1$ MISs. For $s\ge4$, Theorem~\ref{thm:MM} provides a disconnected graph $R$ with $\mis(R)=g(s)$, taking $R=2K_2$ when $s=4$. Since $R$ has no universal vertex, $\mis(R)-\ind_{\{u\in U(R)\}}=g(s)$ for every $u\in V(R)$. The upper bound follows from $\mis(R)\le g(s)$.
\end{proof}

Fix a defining partition of $G\in\B_{s,q,t}$ with support $S$ and kernel $J$, and write $R=G[S]$. With $a_j(u)$ and $P$ as defined in Section~\ref{sec:intro}, define
\[
\Phi(J)=|\{I\in\M(J):I\cap S\ne\varnothing\}|,\qquad
\psi(u)=\prod_{\{j:a_j(u)>0\}}a_j(u)\quad(u\in P).
\]

\begin{lemma}\label{lem:kernel}
Let $G\in\B_{s,q,t}$, and let $J$ be the kernel of the chosen partition. Then
\begin{equation}\label{eq:countkernel}
\mis(G)=4^q3^t+\Phi(J).
\end{equation}
Moreover, for fixed $s,q,t$,
\begin{equation}\label{eq:kernelmax}
\max_{G\in\B_{s,q,t}}\Phi(J)=
\begin{cases}
g(s),&q=0,\\
4^q+\gamma_s,&q=1,2.
\end{cases}
\end{equation}
For the chosen partition, $\Phi(J)$ attains the maximum in~\eqref{eq:kernelmax} if and only if the corresponding condition in Theorem~\ref{thm:equality} holds.
\end{lemma}
\begin{proof}
An MIS avoiding $S$ chooses exactly one vertex from each block $C_j$ and $T_i$. It cannot choose two vertices from the same clique. If it chooses none from a clique, then that clique is undominated, since the set avoids $S$ and there are no edges between distinct blocks. Conversely, choosing one vertex from each block gives an independent set, since there are no edges between distinct blocks. The chosen vertices dominate their respective blocks, and a chosen triangle vertex dominates $S$ because $t\ge1$. Thus the set is maximal. There are $4^q3^t$ choices. An MIS of $G$ meeting $S$ avoids all triangle blocks and is therefore an MIS of $J$. Conversely, an MIS of $J$ meeting $S$ dominates all triangle blocks and is therefore maximal in $G$. This proves~\eqref{eq:countkernel}.

To compute $\Phi(J)$, let $I\in\M(J)$ meet $S$ and put $X=I\cap S$. If $|X|\ge2$, every vertex in a $K_4$ block has a neighbor in $X$, since it has at most one nonneighbor in $S$. Hence $I=X$. If \(X\) is maximal in \(J\), then it is maximal in \(R\), since any support vertex that could be added to \(X\) in \(R\) could also be added in \(J\). Conversely, if \(X\) is maximal in \(R\), then it dominates \(S\) in \(R\), and \(|X|\ge2\) ensures that it dominates every vertex of \(J-S\). These are therefore precisely the nonsingleton MISs of $R$.

If $X=\{u\}$ and $u\notin P$, then $u$ is adjacent to every vertex of $J-S$. Hence $\{u\}$ is maximal in $J$ if and only if $u\in U(R)$. If $X=\{u\}$ with $u\in P$, then maximality requires choosing exactly one of the $a_j(u)$ nonneighbors of $u$ in each block $C_j$ with $a_j(u)>0$. Blocks with $a_j(u)=0$ are already dominated by $u$. These choices are independent across blocks, giving $\psi(u)$ extensions. The resulting set is independent because the selected block vertices lie in distinct blocks and are nonadjacent to $u$. It dominates each block in which a vertex was selected, and any selected block vertex dominates $S\setminus\{u\}$ because it has at most one nonneighbor in $S$. Thus every such extension is maximal in $J$. The three cases contribute $\mis(R)-|U(R)|$, $|U(R)\setminus P|$, and $\sum_{u\in P}\psi(u)$, respectively. Adding these contributions gives
\begin{equation}\label{eq:kernelidentity}
\Phi(J)=\mis(R)-|U(R)\cap P|+\sum_{u\in P}\psi(u).
\end{equation}

For $q=0$, we have $J=R$ and $\Phi(J)=\mis(R)\le g(s)$, with the stated equality condition. The bound is attained by choosing $R$ to be a Moon--Moser extremal graph on $s$ vertices. Suppose $q\ge1$. If $P=\varnothing$, then $\Phi(J)\le g(s)<4^q+\gamma_s$: for $s\ge4$ use $\gamma_s=g(s)$, and for $s\le3$ use $g(s)=s$ and $\gamma_s=s-1$. Thus equality requires $P\ne\varnothing$.

For any $u\in P$, we have $|U(R)\cap P|\ge\ind_{\{u\in U(R)\}}$. Together with Lemma~\ref{lem:gamma}, this gives
\begin{equation*}
\mis(R)-|U(R)\cap P|\le\mis(R)-\ind_{\{u\in U(R)\}}\le\gamma_s.
\end{equation*}
When $q=1$, the remaining term in~\eqref{eq:kernelidentity} is
$\sum_{u\in P}\psi(u)=\sum_{u\in S}a_1(u)\le4$.
The resulting upper bound is $4+\gamma_s$, with equality precisely when this sum is four and $\mis(R)-|U(R)\cap P|=\gamma_s$. These are the conditions in Theorem~\ref{thm:equality}(ii).

When $q=2$, write $a_u=a_1(u)$ and $b_u=a_2(u)$, so $0\le a_u,b_u\le4$ for every $u\in S$. If $a_ub_u=0$ and $u\in P$, then $\psi(u)=a_u+b_u<2(a_u+b_u)$. If both entries are positive, then $\psi(u)=a_ub_u\le2(a_u+b_u)$, by averaging $a_ub_u\le4a_u$ and $a_ub_u\le4b_u$; equality holds only at $a_u=b_u=4$. Consequently,
\begin{equation*}
\sum_{u\in P}\psi(u)\le2\sum_{u\in S}(a_u+b_u)\le16.
\end{equation*}
If $\sum_{u\in P}\psi(u)=16$, then both inequalities in the preceding display are equalities, so $a_u=b_u=4$ for every $u\in P$. The sum is positive, so $P$ is nonempty. Since $\sum_{u\in S}a_u\le4$ and $\sum_{u\in S}b_u\le4$, $P$ contains exactly one vertex $u$. Thus every vertex of both blocks has $u$ as its unique nonneighbor in $S$. Now $P=\{u\}$, so the equality $\mis(R)-|U(R)\cap P|=\gamma_s$ is exactly the support condition in Theorem~\ref{thm:equality}(iii). Conversely, these conditions give equality throughout. This proves the bound $16+\gamma_s$ and all its equality cases.

To attain the bounds for $q=1,2$, choose $(R,u)$ attaining the maximum in Lemma~\ref{lem:gamma}, and make every $C_j$ vertex adjacent to $S\setminus\{u\}$ and nonadjacent to $u$. Then $P=\{u\}$ and $\psi(u)=4^q$, giving $\Phi(J)=\gamma_s+4^q$ by~\eqref{eq:kernelidentity}.

Formula~\eqref{eq:countkernel} also shows that $\Phi(J)=\mis(G)-4^q3^t$ is independent of the choice of defining partition. Hence the corresponding condition in Theorem~\ref{thm:equality} holds for one defining partition if and only if it holds for every defining partition with the same parameters.
\end{proof}

\begin{remark}\label{rem:equality}
Let \(G\) be an extremal graph. When \(q=1\), the four missing edges between $C_1$ and $S$ may be incident with more than one support vertex. When \(q=2\), the missing edges between $C_1\cup C_2$ and $S$ must all be incident with a single vertex of $S$. When $s\ge4$, the support condition in Theorem~\ref{thm:equality}(ii) is equivalent to $\mis(R)=g(s)$ and $U(R)\cap P=\varnothing$; in Theorem~\ref{thm:equality}(iii), it is equivalent to $\mis(R)=g(s)$ and $u\notin U(R)$. The support graphs attaining $g(s)$ are listed in Theorem~\ref{thm:MM}.
\end{remark}

\begin{proof}[Proof of Theorems~\ref{thm:main} and~\ref{thm:equality}]
Fix $k\ge3$, put $s=k-2$, and take $n$ large enough that $t\ge\max\{1,t_0(s)\}$ and $3t\ge k$. Every graph in $\B_{s,q,t}$ then has minimum degree at least $k$, by the degree bounds following Definition~\ref{def:B}. By Lemma~\ref{lem:kernel}, the maximum of $\mis(G)$ over $G\in\B_{s,q,t}$ is the right-hand side of~\eqref{eq:main}.

For a graph outside $\B_{s,q,t}$ with $\delta(G)\ge k$, Lemma~\ref{lem:gap} gives $\mis(G)\le4^q3^t-c_s3^t<4^q3^t$. Since the maximum within $\B_{s,q,t}$ exceeds $4^q3^t$, no admissible graph outside this family is extremal. By~\eqref{eq:countkernel}, a graph in $\B_{s,q,t}$ is extremal if and only if its kernel count $\Phi(J)$ attains the maximum in~\eqref{eq:kernelmax}. By Lemma~\ref{lem:kernel}, this is equivalent to admitting a defining partition that satisfies the corresponding condition in Theorem~\ref{thm:equality}(i)--(iii). Since $t\ge1$, every extremal graph contains triangle vertices of degree exactly $s+2=k$, proving the last assertion of Theorem~\ref{thm:main}.
\end{proof}

\section*{Acknowledgments}
This research was supported by the National Key R\&D Program of China (Grant No. 2023YFA1010202), the National Natural Science Foundation of China (Grant No. 12401448), Fundamental Research Funds for the Central Universities
WK0010000073, and the Natural Science Foundation of Fujian Province (Grant No. 2024J08030).

\section*{Declaration on the use of generative AI}
The authors used ChatGPT to assist in discussing proof strategies, checking proofs, and
improving the exposition. The authors take full responsibility for the mathematical arguments, results,
and conclusions, all of which were carefully reviewed and verified by them.

\Needspace{8\baselineskip}

\end{document}